%% file: paper.tex
\documentclass[journal]{IEEEtran}
				
\usepackage[cmex10]{amsmath}
\usepackage[T1]{fontenc}
\usepackage{letterspace}
\usepackage{amssymb}
\usepackage{amsfonts}
\usepackage{mathtools}
\usepackage{booktabs}
\usepackage{psfrag}
\usepackage{array}
\usepackage{multirow} 
\usepackage{fixltx2e}
\usepackage{cite}

\usepackage{graphicx}
\usepackage[printonlyused,withpage]{acronym}
\usepackage{commath}
\usepackage[font=footnotesize,caption=false]{subfig}
\usepackage{pstool}
\usepackage{tikz,pgfplots}
\usepackage[mediumspace,mediumqspace,squaren,binary]{SIunits}
\usetikzlibrary{shapes,fit,arrows,decorations.markings,matrix,plotmarks,positioning,calc}

\usepackage{ragged2e}
\usepackage{fourier}
\usepackage{setspace}
\usetikzlibrary{chains,shapes.arrows}

\pgfplotsset{
	compat=newest,
	/pgfplots/ylabel absolute/.style={%
		/pgfplots/every axis y label/.style={at={(0,0.5)},xshift=-5pt,rotate=90},
		/pgfplots/every y tick scale label/.style={at={(0,1)},above right,inner sep=0pt,yshift=0.3em}
	}
}%

\newcommand{\ie}{\emph{i.e.},~}
\newcommand{\eg}{\emph{e.g.},}

\newcommand{\mc}{\mathcal}
\newcommand{\ts}[1]{{\textnormal{#1}}}
\newcommand{\pG}{p^{\textnormal{G}}}%
\newcommand{\pW}{p^{\textnormal{W}}}%
\newcommand{\PW}{P^{\textnormal{W}}}%

\newcommand{\PD}{P^{\textnormal{D}}}%
\newcommand{\pD}{p^{\textnormal{D}}}%
\newcommand{\pL}[1]{p^{\textnormal{L}}_{#1}}%
\newcommand{\SG}{\mathcal{G}}%
\newcommand{\SD}{\mathcal{D}}%
\newcommand{\SL}{\mathcal{L}}%
\usetikzlibrary{shapes.geometric,arrows}
\usetikzlibrary{intersections}
\newcommand{\etc}{\emph{etc}}

\title{A Multiperiod OPF Model Under Renewable Generation Uncertainty and Demand Side Flexibility}%
\author{W.~A.~Bukhsh,~\IEEEmembership{Member, IEEE},~C.~Zhang,~\IEEEmembership{Member, IEEE},~P. Pinson,~\IEEEmembership{Senior Member, IEEE}%
\thanks{This work is supported by the research grant from Danish-Chinese project PROAIN.}%

\thanks{W.~A.~Bukhsh, C.~Zhang, and P.~Pinson are with the Department of Electrical Engineering, the Technical University of Denmark, Building 325, 2800 Kgs. Lyngby, Denmark.

(e-mail:~\texttt{\{bukhsh, chzh, ppin\}@elektro.dtu.dk})
}}%

\begin{document}
\maketitle
\begin{abstract}

Renewable energy sources such as wind and solar have received much attention in recent years and large amount of renewable generation is being integrated to the electricity networks. A fundamental challenge in power system operation is to handle the intermittent nature of the renewable generation. In this paper we present a stochastic programming approach to solve a multiperiod optimal power flow problem under renewable generation uncertainty. The proposed approach consists of two stages. In the first stage operating points for conventional power plants are determined. Second stage realizes the generation from renewable resources and optimally accommodates it by relying on demand-side flexibility. The benefits from its application are demonstrated and discussed on a 4-bus and a 39-bus systems. Numerical results show that with limited flexibility on the demand-side substantial benefits in terms of potential additional re-dispatch costs can be achieved. The scaling properties of the approach are finally analysed based on standard IEEE test cases upto 300 buses, allowing to underlined its computational efficiency.
\end{abstract}


\section*{Nomenclature}
{\setstretch{1.35}

\subsection*{Sets}

\begin{description}[\IEEEsetlabelwidth{$P_g^{\ts{G}}, P_g^{\ts{G}+}$}\IEEEusemathlabelsep]

\item[$\mc{B}$] Buses, indexed by $b$.
\item[$\mc{L}$] Lines (edges), indexed by $l$.
\item[$\mc{G}$] Generators, indexed by $g$.
\item[$\mc{W}$] Renewable generators, indexed by $w$.
\item[$\mc{D}$] Loads, indexed by $d$.
\item[$\mc{D}_0$] Flexible loads, $\mc{D}_0 \subseteq \mc{D}$.
\item[$\mc{B}_l$] Buses connected by line $l$.
\item[$\mc{L}_b$] Lines connected to bus $b$.
\item[$\mc{G}_b$] Generators located at bus $b$.
\item[$\mc{D}_b$] Loads located at bus $b$.
\item[$\mc{S}$] Scenarios, indexed by $s$.
\item[$\mc{T}$] Discrete set of time intervals, indexed by $t$.
\end{description}
\subsection*{Parameters}
\begin{description}[\IEEEsetlabelwidth{$P_g^{\ts{G}},P_g^{\ts{G}},
 P_g^{\ts{G}+}$}\IEEEusemathlabelsep]

\item[$b_l$] Susceptance of line $l$.
\item[$\tau_l$] Off-nominal tap ratio of line $l$ (if transformer).
\item[$P_{g}^{\ts{G}-}, P_{g}^{\ts{G}+}$] Min., max. real power outputs of conventional generator $g$.
\item[$\PD_{d,t}$] Real power demands of load $d$.
\item[$f_{g,t}(p^\ts{G}_{g,t})$] Generation cost function for generator $g$.
\item[$\PW_{w,s,t}$] Renewable generation under scenario s from generator w.
\item[$\lambda_{w,s}$] Probability of scenario $s$.
\item[$C^\ts{W}_{w,t}$] Cost of renewable generation spillage.
\item[$F_{d,t}^{-},F_{d,t}^{+}$] Min., max load flexibility of demand at bus $d$.
\item[$\Delta P_{g,t}^{-},\Delta P_{g,t}^{+}$] Min., max change in operating point of generator $g$ during time period $[t,t+1]$.
\end{description}

\subsection*{Variables}
\begin{description}[\IEEEsetlabelwidth{$p_l^{ij}, q_l^{ij}$}\IEEEusemathlabelsep]
\item[$p_{g,t}^\ts{G}$] Real power output of generator $g$.
\item[$p_{w,s,t}^\ts{W}$] Real power output of renewable generator $w$.
\item[$\theta_{b,s,t}$] Voltage phase angle at bus $b$.
\item[$\pL{l,s,t}$] Real power injection at bus $b$
  into line $l$ (which connects buses $b$ and $b'$).
\item[$p_{d,s,t}^\ts{D}$] Real power supplied at bus $d$.
\item[$\alpha_{d,s,t}$] Proportion of load supplied at bus $d$.
\end{description}
}

\section{Introduction}

\IEEEPARstart{E}{lectricity} networks around the world are evolving at a rapid pace. This change is happening because of the increased emphasis on clean and renewable energy sources. Large-scale renewable energy sources (RES) are encouraged by different incentive schemes, in order to support energetic independence and mitigate issues related to climate change. Many countries are investing substantial resources in planning and expanding current infrastructure to cope with RES integration. Wind power generation is the most widely used source of renewable energy and it is been integrated in many power systems around the world \cite{Hawkins,Hatziargyriou}, while solar power is catching up at a rapid pace.

The non-dispatchable nature of wind power introduces additional costs stemming from the management of intermittency \cite{Jabr, Joseph}. Extra reserves need to be obtained, at an additional cost, in order to hedge the uncertainty from the partly predictable generation from wind farms. Despite the advancements in forecasting methodologies and tools, hour-ahead forecast errors for a single wind plant may be as high as 10\%-15\% of its actual output on average \cite{NREL}. 

In contrast, demand at the transmission level has a large base component that can be predicted accurately. In power systems optimization problems electricity demand typically is modelled as inelastic. However in reality a substantial amount of electricity demand is elastic \cite{Baldick}. Electric loads like PEV charging, district heating, HVAC systems are some examples of flexible demands and constitute considerable percentage of the total demand \emph{e.g.}, more than one third of the US residential demand is flexible \cite{USReport}. Based on historic data distribution companies have good idea about the amount of flexibility in demand for a given time window. Majority of these demands are \emph{deferrable} meaning that part of the demand can be shifted in time while respecting deadlines and rate constraints \cite{Baldick}. Demand side management is an active area of research in power systems and there are many challenges with respect to operation and application of demand control \cite{Koutsopoulos}. However demand side flexibility contributes to improve reliability of power system and also reduce the curtailment of renewable sources \cite{Martin}. 

Many distribution system operators implement some sort of demand response (DR) programs. DR can be in the form of interruptible load contracts (ILC) or voluntarily reduction of demand by customers in response to high prices of electricity. DR programs not only provide ancillary services to the operation of power systems but also acts to optimize the generation from renewable sources. Authors in \cite{Tahersima} propose a DR approach and test it on a realistic test case. A stochastic unit commitment model with deferrable demand is presented in \cite{Papavasiliou2}. Authors show that the curtailment of renewable energy can be minimized by initiating DR. 

Traditional power system operation is based on deterministic security-constrained commitment and dispatch models \cite{Allen}. In order to ensure security of supply these models use very conservative forecasts of wind power generation \cite{Hawkins} and as a result of the conservative operation large amount of wind is curtailed \cite{NREL2}. With the increase in wind power penetration in the existing system it is becoming a big challenge to optimally utilize these resources.

The most important decision for a power system operator in short time scales is to determine the operating point of conventional generators (coal, nuclear, \etc) \cite{Phan}. It is very difficult to command an entirely different set point to these generators in short time scales; however small adjustments can be made. Power system operators buy reserve capacities from the fast response units and bring them online if there is a large deviation in the demand. However we note that generation from wind can change considerably in small time scales. Locational marginal prices (LMPs) go up if reserve capacities of generators are utilized. Similarly if more than expected power is generated from wind farms, it has to be spilled in order to keep demand-generation balance; unfortunately resulting in wastage of cheap and clean energy.


The conventional optimal power flow (OPF) problem \cite{carpentier} consists in determining the operating point of generators which minimizes the cost of generation and respects the network and physical constraints. The OPF problem hence provides the dispatch for the next time period, which is usually one-hour ahead. For smaller time scales (typically 5 min) any demand/generation mismatch is alleviated by automatic generation control (AGC) \cite{Jiaqi}. 

The OPF problem has been extended to account for the variable and partly-predictable nature of wind power generation in \emph{e.g.,} \cite{Jabr, Zhang, Vrakopoulou, Entriken, Saunders}. These papers capture the intermittent nature of wind power generation using different probabilistic techniques and determine a robust operating point for the generating units. With stronger focus on the demand side, the authors in \cite{Vrakopoulou} consider demand-side participation as well as uncertainty in demand bids. Finally, in a spirit similar to our proposal, the authors in \cite{Phan} extended the OPF problem to a two-stage stochastic optimization problem, where the decision problem is then to find the steady-state operating point for large generation units in the first stage, while scheduling fast-response generation at the second stage, based on a set of scenarios for renewable energy generation. Demand is there assumed to be deterministic and the problem is not time coupled. This means the optimal operating points are independent of the temporal characteristics of the system.

In view of such limitations in the literature, we propose here to place emphasis on the Multiperiod OPF (MPOPF) problem, which comprises of the time-coupled version of OPF problem. The objective is thus to minimize the cost of generation over the given time horizon while satisfying network constraints and ramp-rate constraints. As an example of recent developments related to the MPOPF, the authors in \cite{Rabiee} consider a stochastic MPOPF model and model the offshore renewable generation with HVDC connections. Uncertainties in wind power generation is considered using a scenario-based approach and demand is assumed to be deterministic.

Stochastic programming approaches \cite{Birge} provide a suitable framework to accommodate the uncertainty in power generation from RES. In this paper we present two stage stochastic program. We consider the flexibility in demands and different scenarios of power production from RES. We focus our attention to intra-hour time scale because of the three reasons. First reason is that the forecasts of renewable generation (mainly wind) are somewhat reasonably accurate in hour ahead period as compared to the day ahead period. Secondly though the amount of energy cleared in short time scales is small but the value of energy is very high. And thirdly given the increasing focus on RES, the penetration from renewable sources would increase and hence the amount of energy cleared in hour-ahead operations. 

The decision problem in our stochastic programming approach is to find the operating point of conventional generators while taking into account the uncertainty in the power generation from RES. Demand flexibilities are considered and optimization decides the operating point of generators, utilization of flexibilities while minimizing total cost of generation. 

Contribution of this paper is to present a framework that can be used to optimally utilize the generation from intermittent sources. Taking flexibilities from demand side and considering possible scenarios of generations from wind power, the proposed approach optimally shifts the demand in the given time horizon. It can also be used as a tool to project future LMPs given demand side flexibilities. The projected prices are useful information for distribution companies, and they can use this information to plan their demand response strategies \cite{Haiwang}. We provide wind scenarios and network data of all the numerical results presented in this paper in an online archive at \cite{DataOnline}.

This paper is arranged as follows. Section \ref{problemformulation} gives the formulation of the problem. Numerical results are given in section \ref{results}. We give conclusions and future research directions in section \ref{conclusion}.

\section{Problem Formulation}\label{problemformulation}

We propose a two-stage stochastic programming formulation of multiperiod optimal power flow problem. In the first stage, decision is made about the dispatch from conventional generators and these decisions remain fixed in the second stage of the problem. The second stage realizes the generation from renewable sources. Any resulting supply-demand mismatch is alleviated by the demand response from flexible demands. It is important to note that demand response can be replaced by the high marginal-cost generators which can be tapped in short term, or the virtual generation resources. However for the sake of clarity of presentation we only consider flexibility in demands. The expensive generators are modelled as high cost generation units in the data, and hence they are minimized as part of the problem.

Consider a power network with set of buses $\mc{B}$. Let $\mc{W}$ denotes the set of renewable generators in the network. Since the real power generation from renewable generators is uncertain, let $\mc{S}$ be the set of real power generation scenarios of these generators. We assume zero marginal price of the generation from renewable generators. Let $\mc{G}$ be the set of conventional power plants. Let $\mc{T}:=\{1,2,\cdots,T\}$ be the set of give time horizon. Following we give constraints and objective function of our two stage stochastic multiperiod optimal power flow problem.

\subsection{Power flow}

Let $\pG_{g,t}$ be the real power generation from the conventional generator $g$ in the time interval $t$. Let $\pW_{w,s,t}$ be the real power generation from the renewable generator $w$ in the time period $t$ in case the scenario $s$ is realized. The power balance equations are given as, $\forall b \in \mc{B}, s \in \mc{S}, t \in \mc{T}$: 
\begin{equation}\label{KCLp}
\sum_{g \in \SG_b} \pG_{g,t}+\sum_{w \in \mc{W}_b} \pW_{w,s,t} = \sum_{d\in \SD_b} \pD_{d,s,t} + \sum_{l\in \SL_b}\pL{l,s,t},
\end{equation}

\noindent
where $\pL{l,s,t}$ is the flow of real power in the line $l$, in the time interval $t$ given scenario $s$ is realized. The power balance equation is given as, $\forall l \in \mc{L}, s \in \mc{S}, t \in \mc{T}$:
\begin{equation}\label{KVLp}
\pL{l,s,t}=  -\frac{b_l}{\tau_l}\left(\theta_{b,s,t}-\theta_{b',s,t}\right),
\end{equation}

\noindent
where $b$ and $b'$ are two ends of the line $l$. Note that we consider the DC model of line flow \cite{jizhong}. This model ignores line losses and reactive power. We have made this assumption in order to keep the formulation linear.

\subsection{Load Model}

Let $\mc{D}$ denotes the set of real power demands and we assume that a distribution network is attached to each bus $d \in \mc{D}$. The demand at distribution network is aggregated and is denoted by $\PD_{d,t}$. We assume that each distribution company at the demand bus $d$ know about the flexibility of their demand during the time interval $t$. This flexibility can either come from distribution company's direct control over some demands or from its DR programs.

Let $\alpha_{d,s,t}$ be the proportion of load supplied to the bus $d$ at the time interval $t$ if the scenario $s$ is realized. Let $[F_{d,t}^{-},F_{d,t}^{+}]$ be the flexibility interval of the demand at bus $d$ during time period $t$. The flexibility in demand can be modelled in different ways. If the demand at the distribution network $d$ is not flexible then $F_{d,t}^{-}=F_{d,t}^{+}=1$ is used. If demand at bus $d$ is flexible then it is placed in the set $\mc{D}_0 \subseteq \mc{D}$. 

The load model is given by following set of constraints:

\begin{subequations}\label{eq:IPload}%
\begin{gather}
\pD_{d,s,t}=\alpha_{d,s,t}\PD_{d,t},\label{eq:loadshedP}\\
0 \le F_{d,t}^{-} \le \alpha_{d,s,t} \le F_{d,t}^{+},\label{eq:loadshedQ}\\
\alpha_{d,s,t}=1, \forall d \in \mathcal{D} \setminus \mathcal{D}_0.
\end{gather}
\end{subequations}

$(1-F_{d,t}^{-})$ is the proportion of demand $d$ which is flexible in the time interval $t$, and $(F_{d,t}^{+}-1)$ is the amount of load that can be increased in the time interval $t$.

 Cost of generation is monotonically increasing function of real power generation. If demands are flexible but not conserved over the given time interval then optimal solution is to reduce the demands. Therefore it is reasonable to consider shifting the demand over a given time period. Distribution company will give its flexibility for each time interval, and the optimization model will decide how to optimally shift the demand. The following constraints ensure that the total demand is met at the end of the time horizon, $~ \forall d \in D_0$:

\begin{equation}\label{sumload}
\sum_{t=1}^T \pD_{d,t} = \sum_{t=1}^T\PD_{d,t}.
\end{equation}

Optimization model would decide the amount of demand to be consumed in each time interval. Note that we assume that there is enough power to support a task which requires more than one time interval to finish. This assumption is justifiable because of the lower bound on the value of $\alpha_{d,s,t}$. Otherwise it is possible to impose a constraint coupled in time. We have assumed that flexibility can be utilized in any way across the time interval. In practice the flexibilities depend on the type of demands \eg~ some demands might need up and down times, and charging/discharging rates. All these technical details can be modelled using linear constraints. However technical details and discussion on this subject is out of the scope of this paper.

\subsection{Operating constraints}

The generation from conventional generators is bounded by the following inequality constraints:

\begin{equation}\label{Pbounds}
P^{\text{G-}}_{g} \le p_{g,t} \le P^{\text{G+}}_{g},
\end{equation}

\noindent
where $P^{\text{G-}}_{g}, P^{\text{G+}}_{g}$ are the lower and upper bounds on the generation output of generator $g$, respectively.

In short time scales it is not be possible for a conventional generator $g$ to considerably deviate from current operating point \cite{Phan}. Therefore we limit the amount of change in generation depending on the ramp rate of individual generators. The constraints are given as:

\begin{subequations}\label{ramp}
\begin{gather}
\Delta P^{-}_{g,t} \le \pG_{g,t+1}-\pG_{g,t},\\
\pG_{g,t+1}-\pG_{g,t} \le \Delta P^{+}_{g,t} .
\end{gather}
\end{subequations}

\subsection{Scenarios of renewable energy generation}
Forecasting of renewable energy generation is a very active area of research, especially for wind and solar energy applications. While forecasts were traditionally provided in the form of single-valued trajectory informing of expected generation for every lead time and location of interest, individually, emphasis is now placed on probabilistic forecasts in various forms \cite{pinson2013}. For decision problems where the space-time dependence structure of the uncertainty is important, forecasts should optimally take the form of space-time trajectories. For example recently, the benefits from employing space-time trajectories in a network-constrained unit commitment problems were demonstrated and discussed \cite{Papavasiliou2}.

In the present case, scenarios of wind power generation are used as input to the stochastic programming approach to solving the multiperiod optimal power flow problem. The exact setup, data and methods of \cite{pinson2013} are employed. In short, the approach relies on nonparametric forecasts for the marginal predictive densities, and on a Gaussian-based copula for the interdependence structure, tracked in an exponential smoothing framework. A sample of 100 space-time scenarios originally issued for 15 control areas in Denmark are used. If others were to aim at reproducing presented results or use these scenarios as input to other stochastic optimization problems, these wind scenarios are made publicly available in an online archive at \cite{DataOnline}.

\subsection{Objective function}
The objective of our optimization is to minimize the cost of generation and optimally utilize the renewable generation. We assume zero marginal price for the renewable generation resources. If such assumption is made in the usual multiperiod OPF problem then it optimizes the renewable resources. However as the demand is fixed, usual formulation of multiperiod OPF will not try to optimally use the flexibility of demand depending on the generation from renewable sources.

Let $\lambda_{w,s}$ be the probability of scenario $s$ for the renewable generator $w$. Also let $C^\ts{W}_{w,t}$ be the cost of wind spillage from generator $g$, in the time interval $t$ respectively. Our objective is to minimize the cost of generation from conventional generators, and optimally utilize the generation from renewable resources while initiating demand response from the distribution system operators. Overall the objective function is to minimize the following over the given time horizon:

\begin{align}\label{obj}
z=\sum_{g \in \SG} f(\pG_{g,t})+\sum_{s \in \mc{S}} \lambda_{w,s}&\left(\sum_{w \in \mc{W}} C^\ts{W}_{w,t}\left(\PW_{w,s,t}-\pW_{w,s,t}\right)\right). 
\end{align}

It is possible to have a cost term in the objective function for shifting demand. However we have assumed that the demand can be shifted freely in the given time horizon. This assumption is based on the understanding that distribution companies will benefit with reduced real time prices if they provide flexibilities in their demands.

\subsection{Overall formulation}
Overall formulation of the multiperiod optimal power flow problem is given as follows:
\begin{subequations}\label{exactOPFpolar}
\begin{gather}
\min \sum_{t \in \mc{T}}z\left(\pG_{g,t},\pW_{w,s,t}\right)\\
\text{subject to}~~~~~~~~~~~~~~~~~~~~~~~~~~~~~~~~~~~~~~~~~~~~~~~~~~~~~~~~~~~~~~~~~~~~~~~~~~~~~~~~~~~~~~~\nonumber\\
(\ref{KCLp}-\ref{ramp}), \label{overall}\\
\theta_{b_0,s,t} = 0,\label{slackbus} \\
0 \le \pW_{w,s,t} \le \PW_{w,s,t},\label{windSpillage}
\end{gather}
\end{subequations}

\noindent
where constraints (\ref{overall}) gives the load model, power balance and power flow equations, bounds of real power generations and ramp rates, respectively. Constraint (\ref{slackbus}) is the slack bus constraint, and (\ref{windSpillage}) is the bounds on the generation from renewable generation.

The overall problem is then, depending on the objective function $f(\pG_{g,t})$ is linear or quadratic program (LP or QP). We use CPLEX 12.06 \cite{CPLEX} called from an AMPL \cite{AMPL} model to solve the problem.

\section{Numerical example}\label{results}

\subsection{An illustrative example: 4 Bus Case}
\begin{figure}
\centering
\input{WB4_windfarmFig.tikz}
\caption{4 Bus Network.}
\label{windfarm}
\end{figure}
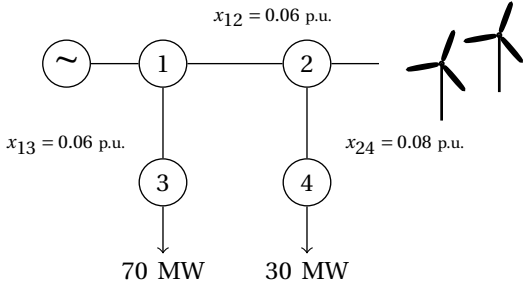

We start with a small 4 bus network as shown in Fig.~\ref{windfarm}. This network consists of one generator and one wind farm. The total load of the network is 100 MW. Complete data of this network is available online at \cite{DataOnline}. 

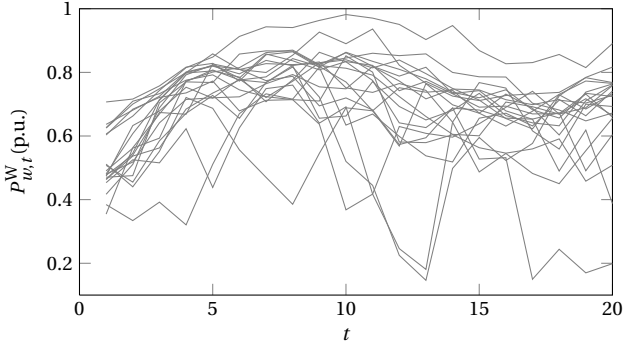
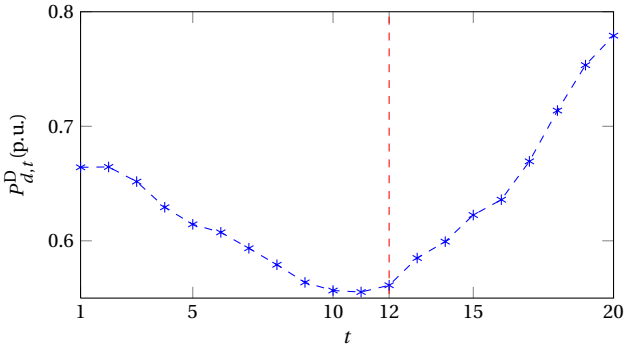
\begin{figure}
\centering\footnotesize
\subfloat[20 Wind Generation Scenarios.]{\input{windscenF.tikz}}\\
\subfloat[Load Profile.]{\input{loadcurveF.tikz}}
\caption{Wind scenarios and demand profile for the 4 bus network.}
\label{windscen}
\end{figure}

We assume that the time horizon consists of twenty time periods \emph{i.e.} $\mc{T}=\{1,2,\cdots,20\}$. We assume 20 different scenarios for wind power generation at bus 1 as shown in Fig.~\ref{windscen}(a). 

We assume zero marginal cost for the wind power. The marginal price of conventional generator at bus 1 is nonzero and quadratic monotonically increasing function of real power generations. We assume the cost of wind spillage to be unity and ramp rate of the generator at bus 1 to be $\pm 10\%$. It is important to note that for feasibility the least value of ramp rate should be equal or greater than the max rate of change in demand during any given time interval. For this test case the maximum change of 6\% occurs between the time periods 18 and 19 (Fig.~\ref{windscen}(b)).

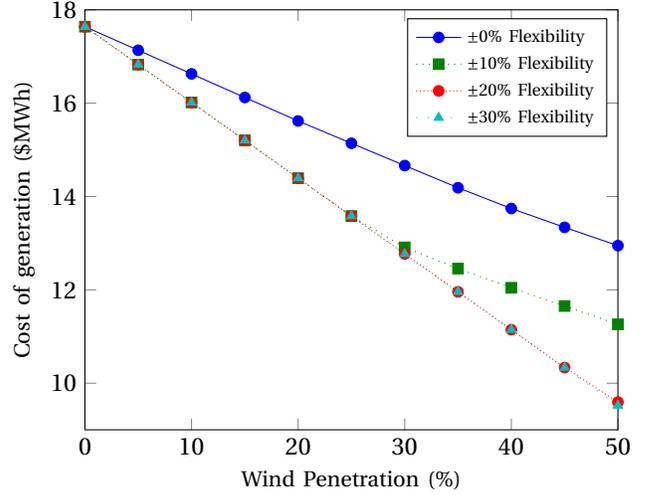
\begin{figure}[t]
\centering
\input{2busgencost_20scen.tikz}
\caption{Generation cost vs wind power penetration for 4 bus network.}
\label{generationcost}
\end{figure}

Fig.~\ref{generationcost} shows the cost of generation as the wind power penetration is increased in the system. We can observe that the cost of generation is monotonically decreasing as the wind power penetration in the system is increased. Also the cost of generation decreases further when the demand is made more flexible. There is no difference in the cost of generation between $\pm 20\%$ and $\pm 30\%$ demand flexibility. This is because ramp rate of the conventional generator is not fast enough to utilize the flexibility of demand. For this example we can say that for given ramp rate of $\pm 10\%$ and wind generation uncertainties, the optimal demand flexibility needed to fully utilize the wind power is $\pm20\%$.

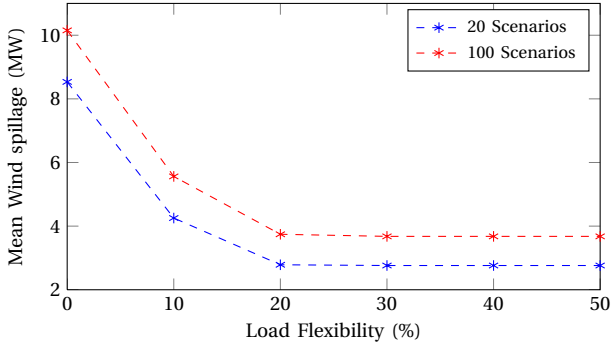
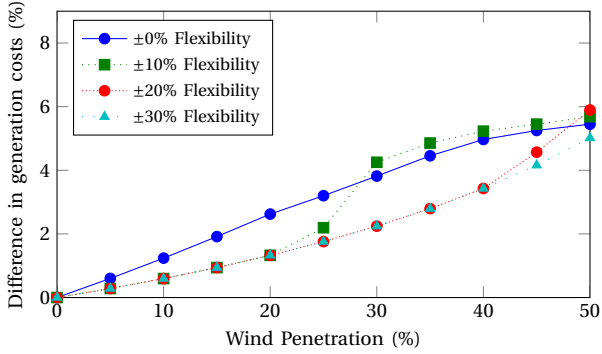
\begin{figure}[t]
\centering\footnotesize
\subfloat[Wind spillage.]{\input{windspill.tikz}}\\
\subfloat[Value of Stochastic Solution (20 scenarios vs 100 Scenarios).]{\input{valueofstochasticsol.tikz}}
\caption{Robustness of the solutions of 4 bus network with respect to uncertainty in the wind power generation.}
\label{valueofstoch}
\end{figure}

Cost of generation depends upon the uncertainty in the wind power generation. If we increase the number of scenarios then the cost of generation would increase. Fig.~\ref{valueofstoch} shows the robustness of solution depending on the number of scenarios. We increase the number of scenarios from 20 to 100 and we can observe in Fig.~\ref{valueofstoch}(a) that the mean wind spillage (for all scenarios and all time periods) is increased. Fig.~\ref{valueofstoch}(b) shows the difference in cost of generation. The difference in cost of generation between 20 and 100 scenarios increases as the wind penetration in the system increase. This is because there is more uncertainty in generation from wind for 100 scenarios as compared to 20 scenarios. However the difference between cost of generation, for given demand flexibilities and penetration levels, is always less than 6\%.

\subsection{39 Bus Case}

\begin{figure}
\includegraphics[width=0.9\columnwidth ]{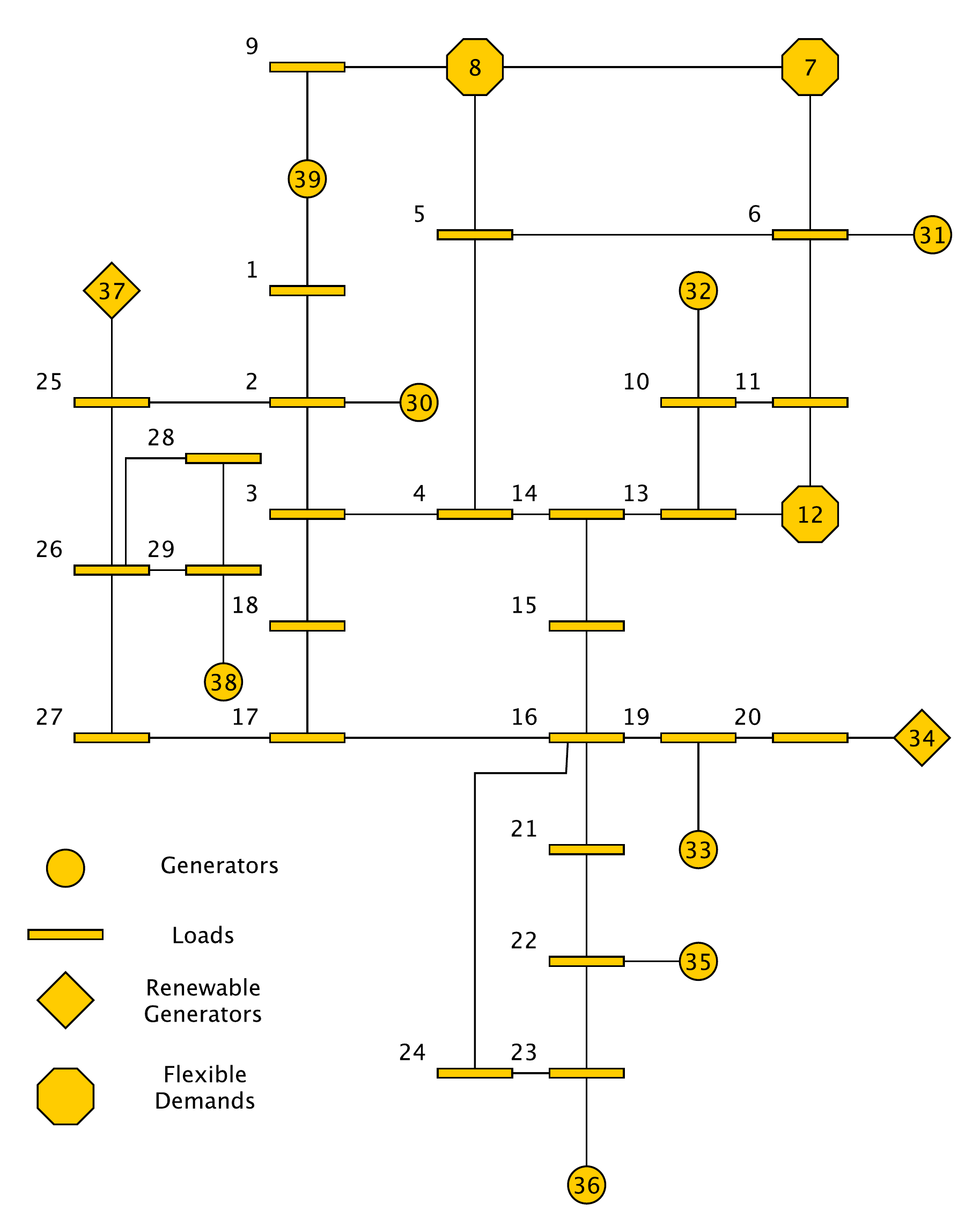}
\caption{Modified 39 bus system}
\label{case39}
\end{figure}

Consider the 39 bus New England test network obtained from \cite{zimmerman}. This test network consists of 39 buses, 10 generators, and 46 transmission lines. We modify the network as follows. We consider there are 8 conventional generators, and there are two renewable sources at buses 34 and 37 respectively. We consider that demands at buses 7, 8 and 12 are flexible demands. The topology of the network is shown in Fig.~\ref{case39}. Default data from \cite{zimmerman} assume same cost data for all generators. We take more realistic cost data from \cite{39busdata} to use in our example. Modified data of this network is available at \cite{DataOnline}.

The total demand in the network is 6254.23 MW. Approximately 27\% of this demand is at the flexible demand buses 4, 8 and 20. Total generation capacity of the network is 7367 MW, and approximately 15\% of the total capacity is from renewable generators at buses 34 and 37. We assume the ramp rate of conventional generators to be $\pm 5\%$. 

\begin{figure}[t]
\centering\footnotesize
\subfloat[System prices.]{\input{LMPS_39bus.tikz}}\\
\subfloat[Expected demand depending on flexibility.]{\input{loadprofile_39bus.tikz}}
\caption{Numerical results for 39 bus system.}
\label{39bus_results}
\end{figure}
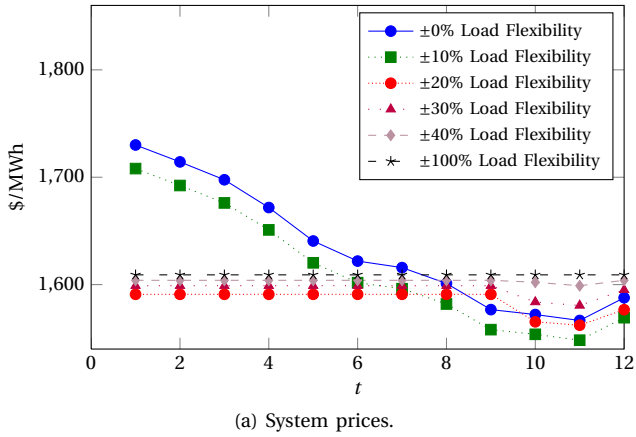
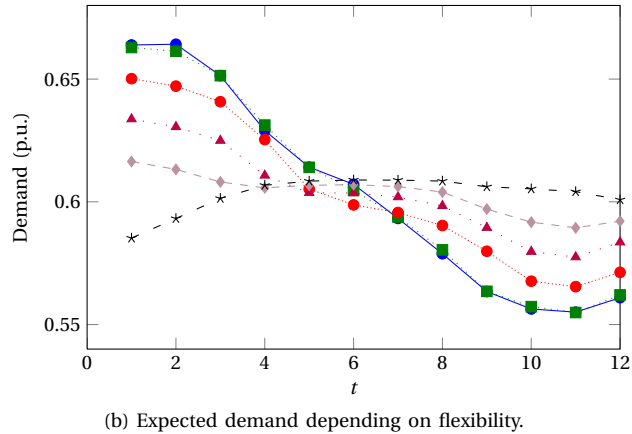

Let $\mc{T}=\{1,2,\cdots,12\}$ (first 12 from Fig.~\ref{windscen}(b)) and consider $100$ independent scenarios for the renewable generators at the buses 34 and 37. Fig.~\ref{39bus_results} shows the result of our model as the flexibility of demand is increased. Line limits were not active at the optimal solution, therefore the locational marginal price at all buses were equal. The solid (blue) line shows the results when demand at buses 7, 8 and 12 is not flexible. In this case the marginal prices follow the behaviour of demand curve \ie prices are high when demand is high and prices decrease with the decrease in demand. If demand is $\pm 10\%$ flexible then the marginal prices are low but this flexibility (coupled with $\pm 5\%$ ramp rate) is not enough to have constant system price. We observed that with $\pm 10\%$ demand flexibility, the cost of generation is decreased by $3.9\%$. Further as the flexibility of demand is decreased, the system price tends toward a constant function. It is interesting to note that the difference in system prices is very small for the demand flexibilities of $40\%$ and $100\%$. This is because of the strict ramp rate constraints, \ie generators can not change their operating point fast enough to utilize the flexibility of demand. However in practice it is not plausible to have $100\%$ flexibility at any demand node. 

Another interesting point to observe is that since we consider the linear model of the system, the results are generally independent of the flexibility \ie the flexibility can come from any node of the network as long as line limits are respected. In practice the transmission system are lossy, so the results would depend on line losses however the effect of line losses is expected to be very small.

\subsection{Larger test cases}
We consider the standard IEEE test networks consisting of 14, 30, 57, 118 and 300 buses from the test archive at \cite{WashingtonTest}. We also consider 9, 24 and 39 bus test cases from \cite{zimmerman}. For all test cases we assumed ramp rate of conventional generators to be $\pm 10\%$, number of scenarios to be 50 and 12 time intervals. We generated large number of scenarios by considering different demand flexibilities and choices of wind generation buses. To keep consistency across all scenarios we considered that for all cases wind power penetration is always less than or equal to 25\%. For all the instances total demand across the time horizon is constrained to be conserved.

Tab.~\ref{largertestnetworks} gives the results of some of the scenarios on 57, 118 and 300 bus networks. Second column in this table gives the set of buses where wind power generation is assumed. Third column gives the percentage of wind power penetration in system. Column four and five gives the set of buses which are flexible and their percentage of load in the system respectively. Second last column gives the assumed flexibility in the set $\mc{D}_0$. Last column shows the improvement in the cost of generation when compared to solving the problem with inflexible loads.

Results in Tab.~\ref{largertestnetworks} shows that considerable savings can be made in the generation cost if demands are flexible. For example consider the 57 bus case with $W=\{7\}$ and $\mc{D}_0=\{12\}$. In this case the load at bus $12$ is approximately 30\% of the total load of the network. The result shows that if the demand at bus $12$ is $\pm 10\%$ flexible that the cost of generation can be improved by 4\%, \ie approximately $3\%$ ($10\%$ of $30\%$) flexibility in demand results in $4\%$ reduction in cost of generation.

\renewcommand{\arraystretch}{1.5}
\begin{table*}[t]
\centering\footnotesize
\caption{Results on larger test networks.}
\label{largertestnetworks}
\begin{tabular}{llclccc}
\toprule
&$\mc{W}$&Wind Penetration (\%)&$\mc{D}_0$&$\frac{\mc{D}_0}{\mc{D}}$(\%)&$\mc{D}_0$ Flexibility ($\pm$\%)&Improvement (\%) \\
\midrule
\multirow{4}{*}{\texttt{case57}}&\multirow{2}{*}{\{3\}}&\multirow{2}{*}{7.1}& \{8\}&12.0&10&1.65\\
		&					  & 	           &   \{12\}  &30.1  &10&4.13\\  \cline{2-7}
&					  \multirow{2}{*}{\{12\}}&      \multirow{2}{*}{20.7}&    \{8\}    &12.0 &10&2.00\\  
		&					  & 	           &   \{9\}  &9.7    &20&3.15\\

\midrule
\multirow{4}{*}{\texttt{case118}}&\multirow{2}{*}{\{10\}}&\multirow{2}{*}{5.5}& \{80, 116\}&7.4&20&1.97\\
		&					  & 	           &   \{54\} &2.6     &10&0.36\\  \cline{2-7}
&					  \multirow{2}{*}{\{69, 89\}}&      \multirow{2}{*}{15.2}&    \{42, 59, 90\}  &12.6   &20&4.25\\  
		&					  & 	           &   \{54\} &2.6     &10&0.45\\  

\midrule
\multirow{4}{*}{\texttt{case300}}&\multirow{2}{*}{\{186, 191\}}&\multirow{2}{*}{10.3}& \{5, 20\}&4.1&20&1.13\\
		&					  & 	           &   \{120, 138, 192\} &11.0     &20&3.03\\  \cline{2-7}
&					  \multirow{2}{*}{\{191, 7003, 7049, 7130\}}&      \multirow{2}{*}{21.9}&    \{10,44\} &1.5    &10&0.24\\  
		&					  & 	           &   \{120, 138, 192\}  &6.2    &20&3.45\\  

\bottomrule
\end{tabular}
\end{table*}

\begin{figure}[t]
\centering
\input{times.tikz}
\caption{Min., mean and max. solution times for solving multiperiod OPF with different demand flexibilities and wind penetration.}
\label{timing}
\end{figure}
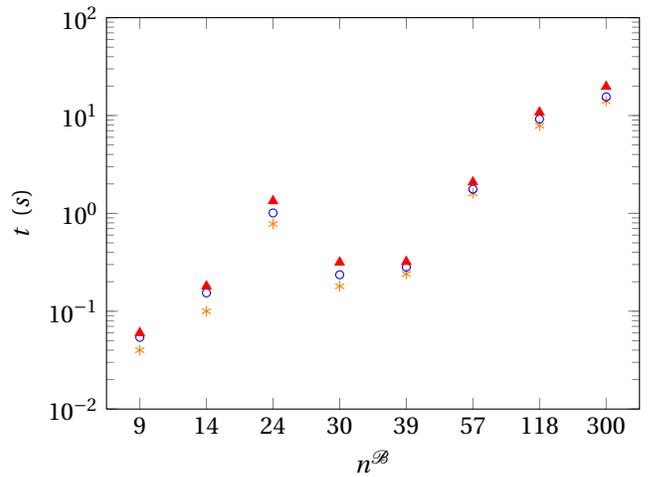

Fig.~\ref{timing} gives the run times on all standard test cases. Problems were solved on a single core 64 bit Linux machine with 8 GiB RAM, using AMPL 11.0 with CPLEX 12.6 to solve LP and QP problems. The results are for large number of scenarios  for wind power penetration (less than 25\%) and demand flexibilities. Fig.~\ref{timing} shows that the solution times scale well with increase in the size of the network. Note that solution times for 24 bus case is higher than 39 and 57 network. This is because of the reason that 24 bus network has more generators than 39 and 57 bus networks and hence the size of the problem is bigger.

\section{Conclusions and Future Research Directions}\label{conclusion}

In this paper we presented a two stage stochastic programming approach to solve multiperiod optimal power flow with flexible demands. We observed that considerable savings in power generation costs can be made if a small proportion of the demand is flexible. The flexibility of the demand can come from any node of the network provided it respects the network constraints. Numerical results show that the uncertain wind power generation can be optimally utilized using flexibility of the demand and hence maximizing the social welfare. Computational times shows the promise of the proposed approach.

Future research will investigate the wider practical aspects of the approach. We would like to extend this to longer time scales by considering unit-commitment as part of the problem. Current research work is looking at modelling this problem with AC power flow equations \ie considering line losses and reactive power flows.
\bibliography{biblio}

\vspace{-0.7cm}

\begin{IEEEbiographynophoto}
{\bf Waqquas A. Bukhsh}
(S'13, M'14) received the B.S. degree in Mathematics from
COMSATS University, Pakistan in 2008 and the Ph.D. degree in Optimization and Operational Research from The University of Edinburgh, UK in 2014.

He is a post doctoral researcher at the Department of Electrical Engineering of the Technical University of Denmark. His research interests are in large scale optimization methods and their application to energy systems.
\end{IEEEbiographynophoto}

\vspace{-0.9cm}

\begin{IEEEbiographynophoto}
{Chunyu Zhang} (M'12) received the B.Eng., M.Sc., and Ph.D. degrees from North China Electric Power University, China, in 2004, 2006 and 2014, respectively, all in electrical engineering. 

From 2006 to 2012, he joined National Power Planning Center and CLP Group Hong Kong as a senior engineer. He is currently pursuing the Ph.D. degree at the Center for Electric Power and Energy, Technical University of Denmark (DTU), Denmark. His research interests include power systems planning and economics, smart grid and electricity markets.
\end{IEEEbiographynophoto}

\vspace{-0.9cm}

\begin{IEEEbiographynophoto}
{Pierre Pinson} (M'11, SM'13) received the M.Sc. degree in Applied Mathematics from the National Institute for Applied Sciences (INSA Toulouse, France) and the Ph.D. degree in Energy from Ecole des Mines de Paris. 

He is a Professor at the Department of Electrical Engineering of the Technical University of Denmark, where he heads the Energy Analytics \& Markets group of the Center for Electric Power and Energy. His research interests include among others forecasting, uncertainty estimation, optimization under uncertainty, decision sciences, and renewable energies. He acts as an Editor for the {\it IEEE Transactions on Power Systems} and for {\it Wind Energy}.
\end{IEEEbiographynophoto}

\vfill

\end{document}

%% file: WB4_windfarmFig.tikz
\newcommand\windturbine[2][]{
  \begin{scope}[shift={(#2)}]
    \draw[#1] (0,0) circle (0.01cm );
    \draw[#1,rotate=70] (0.2,0) ellipse (0.15cm and 0.01cm);
    \draw[#1,rotate=190] (0.2,0) ellipse (0.15cm and 0.01cm);
    \draw[#1,rotate=310] (0.2,0) ellipse (0.15cm and 0.01cm);
  \end{scope}
}

\footnotesize
\tikzstyle{gen} = [draw, circle, minimum size=2em]
\tikzstyle{con} = [draw, circle, minimum size=2em, double]
\tikzstyle{bus} = [draw, rectangle, fill=black, inner sep=0pt, minimum width=4em, minimum height = 0em, very thick]
\tikzstyle{busv} = [draw, rectangle, fill=black, inner sep=0pt, minimum width=0em, minimum height = 4em, very thick]
\tikzset{->-/.style={decoration={
  markings,
  mark=at position #1 with {\arrow{>}}},postaction={decorate}}}
\tikzset{<--/.style={decoration={
  markings,
  mark=at position #1 with {\arrow{<}}},postaction={decorate}}}
\resizebox{2.8in}{1.5in}
{%
\begin{tikzpicture}[]
\node[gen, inner sep=0pt] (k1) {$1$};
\node[gen,right=4em of k1] (k2) {$2$};
\node[gen,below=3em of k1] (k3) {$3$};
\node[gen,below=3em of k2] (k4) {$4$};

\node[gen, left=2em of k1,inner sep=0pt] (G1) {\large$\sim$};
\node[right=2em of k2] (G2) {};

\node[below=2em of k3] (D1) {70 MW};
\node[below=2em of k4] (D2) {30 MW};

\node[above right= 0.5em and 1em of k1] (LL1) {\tiny $x_{12}=0.06$ p.u.};
\node[below left= 2em and 0.5em of k1] (LL2) {\tiny $x_{13}=0.06$ p.u.};
\node[below right= 2em and 0.5em of k2] (LL3) {\tiny $x_{24}=0.08$ p.u.};

%
%
%
\draw[-] (k1) -- (k2);
\draw[-] (k1) -- (k3);
\draw[-] (k2) -- (k4);

\draw[-] (k1) -- (G1);
\draw[-] (k2) -- (G2);

\draw[->] (k3) -- (D1);
\draw[->] (k4) -- (D2);

  \windturbine[draw=black,fill=black, very thick]{2.9,0.0};  
  \draw[thick] (2.9,0.0) -- (2.9,-0.6);
    \windturbine[draw=black,fill=black, very thick]{3.5,0.3};  
      \draw[thick] (3.5,0.3) -- (3.5,-0.3);

\end{tikzpicture}
}

%% file: windscenF.tikz
%
%
%
\definecolor{mycolor1}{rgb}{0.00000,0.75000,0.75000}%
\definecolor{mycolor2}{rgb}{0.75000,0.00000,0.75000}%
\definecolor{mycolor3}{rgb}{0.75000,0.75000,0.00000}%
\begin{tikzpicture}

\begin{axis}[%
width=0.8\columnwidth,
height=1.5in,
scale only axis,
xmin=0,
xmax=20,
ymin=0.1,
ymax=1,
xlabel=$t$,
ylabel=$P^\ts{W}_{w,t}$(p.u.)
]
\addplot [color=gray]
  table[row sep=crcr]{
1	0.384460432271812	\\
2	0.334138265321714	\\
3	0.39211021799672	\\
4	0.320718432808069	\\
5	0.511097833226585	\\
6	0.670195008873801	\\
7	0.732582856766447	\\
8	0.715879043358315	\\
9	0.81648374993383	\\
10	0.863173498473644	\\
11	0.834677137635638	\\
12	0.809602492154417	\\
13	0.779704415873997	\\
14	0.737250632940909	\\
15	0.720228322487059	\\
16	0.74521023455785	\\
17	0.682319241783506	\\
18	0.656484829180529	\\
19	0.734256359587207	\\
20	0.72407359526009	\\
};
\addplot [color=gray]
  table[row sep=crcr]{
1	0.507700265124861	\\
2	0.454994581411481	\\
3	0.584795107848096	\\
4	0.671243011484326	\\
5	0.727537368423861	\\
6	0.655562455470796	\\
7	0.76834863691622	\\
8	0.816961954548812	\\
9	0.722400860955903	\\
10	0.52011102300564	\\
11	0.444935503018041	\\
12	0.224880140759073	\\
13	0.146453658916529	\\
14	0.498119722706043	\\
15	0.595258428570188	\\
16	0.523097379856881	\\
17	0.576699300007956	\\
18	0.726077589122989	\\
19	0.784294396660919	\\
20	0.816551027873552	\\
};
\addplot [color=gray]
  table[row sep=crcr]{
1	0.464001467924984	\\
2	0.545837149872447	\\
3	0.714400614906242	\\
4	0.79785115199653	\\
5	0.804718679908789	\\
6	0.785114409675001	\\
7	0.70783038665257	\\
8	0.776469661762582	\\
9	0.863191232309735	\\
10	0.833182001656301	\\
11	0.774479239469095	\\
12	0.641340054250144	\\
13	0.626471977481607	\\
14	0.687403886844738	\\
15	0.681553459721394	\\
16	0.656027907387665	\\
17	0.653150107884291	\\
18	0.673278470745778	\\
19	0.490773234972851	\\
20	0.604189638184869	\\
};
\addplot [color=gray]
  table[row sep=crcr]{
1	0.476854387845798	\\
2	0.539685807934138	\\
3	0.673901002894006	\\
4	0.668318094370437	\\
5	0.82651631620741	\\
6	0.809137827933815	\\
7	0.799635018185334	\\
8	0.854023534016241	\\
9	0.92659017097977	\\
10	0.890071784167414	\\
11	0.936281888606562	\\
12	0.822665472537787	\\
13	0.765821449581227	\\
14	0.725455842077439	\\
15	0.71141171051656	\\
16	0.727362719356604	\\
17	0.732590527541154	\\
18	0.737030766448878	\\
19	0.783363367970907	\\
20	0.768001548953736	\\
};
\addplot [color=gray]
  table[row sep=crcr]{
1	0.480010190535918	\\
2	0.51808091470278	\\
3	0.641382393207895	\\
4	0.714843169393286	\\
5	0.720152000170949	\\
6	0.72318075693519	\\
7	0.764298474493265	\\
8	0.793214380958502	\\
9	0.754503820232877	\\
10	0.748322057286743	\\
11	0.818448145867091	\\
12	0.772781997692101	\\
13	0.716743790921845	\\
14	0.747518233197063	\\
15	0.744790019574768	\\
16	0.714129446946162	\\
17	0.698494334156966	\\
18	0.714711640992591	\\
19	0.66360674128172	\\
20	0.741312681835563	\\
};
\addplot [color=gray]
  table[row sep=crcr]{
1	0.354536609476722	\\
2	0.563091923283177	\\
3	0.669161280575025	\\
4	0.771166690036516	\\
5	0.807447888972683	\\
6	0.776111427937298	\\
7	0.768020472210043	\\
8	0.821533153487488	\\
9	0.832787657979085	\\
10	0.75147288328809	\\
11	0.737027839766242	\\
12	0.764897387646617	\\
13	0.74345594135004	\\
14	0.643461257995582	\\
15	0.486981317771057	\\
16	0.545111969385291	\\
17	0.482333176548087	\\
18	0.450140018368363	\\
19	0.553687782238319	\\
20	0.655756572044331	\\
};
\addplot [color=gray]
  table[row sep=crcr]{
1	0.603020114597728	\\
2	0.684235601445033	\\
3	0.721244414825588	\\
4	0.735418945657102	\\
5	0.772644927052193	\\
6	0.678846027112364	\\
7	0.712148103016315	\\
8	0.714238653675194	\\
9	0.709938598615483	\\
10	0.852109974916606	\\
11	0.800651349297034	\\
12	0.581162935926627	\\
13	0.579135705163709	\\
14	0.604101460732588	\\
15	0.563243839705589	\\
16	0.546089911980437	\\
17	0.559460530152808	\\
18	0.5909413467728	\\
19	0.458039783823324	\\
20	0.508252063556869	\\
};
\addplot [color=gray]
  table[row sep=crcr]{
1	0.636560532769823	\\
2	0.681173440683894	\\
3	0.724912174757764	\\
4	0.80125739278856	\\
5	0.820793711639521	\\
6	0.808963454687581	\\
7	0.827495629791464	\\
8	0.824394143311526	\\
9	0.828912757525642	\\
10	0.839561793495942	\\
11	0.814169069743811	\\
12	0.838820051378025	\\
13	0.796081412346027	\\
14	0.749488849753376	\\
15	0.713010494951706	\\
16	0.654871162129247	\\
17	0.622464857535683	\\
18	0.729499485149222	\\
19	0.684765899358787	\\
20	0.733457289155256	\\
};
\addplot [color=gray]
  table[row sep=crcr]{
1	0.470492310041842	\\
2	0.475441392844854	\\
3	0.704372636553627	\\
4	0.770115187272298	\\
5	0.768051220555416	\\
6	0.857584709400112	\\
7	0.867575623877164	\\
8	0.86666333876156	\\
9	0.821415256603287	\\
10	0.813039094633973	\\
11	0.791317217091432	\\
12	0.724023536310422	\\
13	0.661394648336115	\\
14	0.605207571156352	\\
15	0.600573977225107	\\
16	0.706721484658384	\\
17	0.691989649769045	\\
18	0.668898954996304	\\
19	0.71706470899175	\\
20	0.664312692456264	\\
};
\addplot [color=gray]
  table[row sep=crcr]{
1	0.511492809732616	\\
2	0.440283092709454	\\
3	0.551684719677384	\\
4	0.721009462413372	\\
5	0.6870239042129	\\
6	0.556367172249504	\\
7	0.470153679455306	\\
8	0.385440955028697	\\
9	0.54528100815822	\\
10	0.688690731559196	\\
11	0.679155400757476	\\
12	0.568450191448247	\\
13	0.764545862686758	\\
14	0.697999868506397	\\
15	0.671704608553886	\\
16	0.631148604982174	\\
17	0.677907271349313	\\
18	0.731008772775424	\\
19	0.772831667029114	\\
20	0.765306955365292	\\
};
\addplot [color=gray]
  table[row sep=crcr]{
1	0.60638890219245	\\
2	0.6626146831542	\\
3	0.743107034452889	\\
4	0.814958765937026	\\
5	0.848104266510175	\\
6	0.912908861134967	\\
7	0.943198565045559	\\
8	0.940802343568089	\\
9	0.963264940085333	\\
10	0.981795368592352	\\
11	0.970537344067776	\\
12	0.950654863249211	\\
13	0.903162037229189	\\
14	0.947121102856235	\\
15	0.868620340505388	\\
16	0.827103335356677	\\
17	0.830768092991951	\\
18	0.856224457630414	\\
19	0.815049856759624	\\
20	0.890773280753169	\\
};
\addplot [color=gray]
  table[row sep=crcr]{
1	0.416816041022155	\\
2	0.510750071898441	\\
3	0.688746606408176	\\
4	0.792234839336501	\\
5	0.823361213428468	\\
6	0.782131454372626	\\
7	0.857831103612313	\\
8	0.836430639589717	\\
9	0.807927007691385	\\
10	0.833836102241968	\\
11	0.821980621095995	\\
12	0.742305688808134	\\
13	0.771709572448347	\\
14	0.732842132050714	\\
15	0.690640456037334	\\
16	0.69390138040507	\\
17	0.669433205743139	\\
18	0.670463135934134	\\
19	0.731855660265879	\\
20	0.659269082564412	\\
};
\addplot [color=gray]
  table[row sep=crcr]{
1	0.490449714620625	\\
2	0.561861083072281	\\
3	0.593205857790327	\\
4	0.77299900665613	\\
5	0.800899761666246	\\
6	0.814821778127037	\\
7	0.86204516137222	\\
8	0.869959811306939	\\
9	0.824121476052643	\\
10	0.861640259427898	\\
11	0.819178672520806	\\
12	0.788329383543514	\\
13	0.816045182638405	\\
14	0.733559564807954	\\
15	0.742588822543514	\\
16	0.679028044287538	\\
17	0.652605550424679	\\
18	0.489830105422221	\\
19	0.618663564813841	\\
20	0.387693464148857	\\
};
\addplot [color=gray]
  table[row sep=crcr]{
1	0.62607518883265	\\
2	0.709355630326193	\\
3	0.732948663500034	\\
4	0.775471455718017	\\
5	0.792235788202152	\\
6	0.745059768305872	\\
7	0.723586200828902	\\
8	0.71424280771285	\\
9	0.637060515573696	\\
10	0.367945408232246	\\
11	0.417361066989685	\\
12	0.629844598934385	\\
13	0.608314656063135	\\
14	0.596334490613096	\\
15	0.638498084711882	\\
16	0.69026687814563	\\
17	0.740590263371003	\\
18	0.668978727322691	\\
19	0.688797431141341	\\
20	0.758632689526005	\\
};
\addplot [color=gray]
  table[row sep=crcr]{
1	0.624891767890761	\\
2	0.700359043393686	\\
3	0.74844114311062	\\
4	0.799157682714638	\\
5	0.858176334641582	\\
6	0.809832975322655	\\
7	0.773981323586879	\\
8	0.787824917430719	\\
9	0.694793204757166	\\
10	0.719651741145495	\\
11	0.68193839387834	\\
12	0.573168846982576	\\
13	0.594569408152306	\\
14	0.688140092717649	\\
15	0.527456185546504	\\
16	0.53196562086283	\\
17	0.149784843207464	\\
18	0.243996968434417	\\
19	0.169851247311245	\\
20	0.198749292122223	\\
};
\addplot [color=gray]
  table[row sep=crcr]{
1	0.453345920245346	\\
2	0.52516898463588	\\
3	0.515607990916382	\\
4	0.622889221828565	\\
5	0.438353477250778	\\
6	0.62551414447557	\\
7	0.748862783882392	\\
8	0.760024997559823	\\
9	0.640460893072119	\\
10	0.6915515140158	\\
11	0.417313608579893	\\
12	0.246017627603893	\\
13	0.180822282294279	\\
14	0.643627326189496	\\
15	0.767888845206917	\\
16	0.750700024499024	\\
17	0.616018907429441	\\
18	0.559003904739398	\\
19	0.700090195482725	\\
20	0.692748834641714	\\
};
\addplot [color=gray]
  table[row sep=crcr]{
1	0.473866968479603	\\
2	0.632762426415793	\\
3	0.678435755115963	\\
4	0.75332281651672	\\
5	0.718846032738388	\\
6	0.779079744333189	\\
7	0.721023926684354	\\
8	0.772188959695576	\\
9	0.829976214845266	\\
10	0.845189336636004	\\
11	0.860171469412401	\\
12	0.852590997641089	\\
13	0.858093819932469	\\
14	0.799798681434136	\\
15	0.786081195519918	\\
16	0.785374298430299	\\
17	0.703105376577638	\\
18	0.681682344439851	\\
19	0.647517889089783	\\
20	0.672212870657989	\\
};
\addplot [color=gray]
  table[row sep=crcr]{
1	0.463208901473906	\\
2	0.536932940233576	\\
3	0.562456714342739	\\
4	0.684585321261392	\\
5	0.73493368294274	\\
6	0.775982730910883	\\
7	0.785816021402743	\\
8	0.843245455073403	\\
9	0.818548072728096	\\
10	0.769763872713641	\\
11	0.679231034693697	\\
12	0.75175093683041	\\
13	0.737661499197194	\\
14	0.743693511946189	\\
15	0.671793108630265	\\
16	0.731673241491226	\\
17	0.687186003679787	\\
18	0.603286254033517	\\
19	0.710622087479164	\\
20	0.760273488083781	\\
};
\addplot [color=gray]
  table[row sep=crcr]{
1	0.480926612017417	\\
2	0.607539680215699	\\
3	0.712530805458232	\\
4	0.799003316107135	\\
5	0.819646490608887	\\
6	0.795299052894718	\\
7	0.837647744434076	\\
8	0.820739853275601	\\
9	0.655836455097739	\\
10	0.662785621966269	\\
11	0.771159168680745	\\
12	0.692811408397315	\\
13	0.650502711628269	\\
14	0.680562763538312	\\
15	0.733950674297549	\\
16	0.664157706625163	\\
17	0.641256599783728	\\
18	0.6320815460347	\\
19	0.738067971748118	\\
20	0.763531887895715	\\
};
\addplot [color=gray]
  table[row sep=crcr]{
1	0.70697671794577	\\
2	0.716141280012251	\\
3	0.757442428600291	\\
4	0.816397933526318	\\
5	0.828288071288077	\\
6	0.772605658179615	\\
7	0.855393373731075	\\
8	0.86460102553616	\\
9	0.815997752128379	\\
10	0.634213395428634	\\
11	0.669105598689948	\\
12	0.598809264556469	\\
13	0.53764726770764	\\
14	0.518361117066259	\\
15	0.655976546706677	\\
16	0.671679024009938	\\
17	0.63158938493826	\\
18	0.712604483353321	\\
19	0.703563712938469	\\
20	0.72756830490798	\\
};
\end{axis}
\end{tikzpicture}%

%% file: loadcurveF.tikz
%
%
\begin{tikzpicture}

\begin{axis}[%
width=0.8\columnwidth,
height=1.5in,
scale only axis,
xmin=1,
xmax=20,
ymin=0.55,
ymax=0.8,
xlabel=$t$,
ylabel=$P^\ts{D}_{d,t}$(p.u.),
xtick={1,5,10,12,15,20},
xticklabels={1,5,10,12,15,20},
]
\addplot [color=blue,dashed,mark=asterisk,mark options={solid},forget plot]
  table[row sep=crcr]{
1	0.6642	\\
2	0.6645	\\
3	0.6517	\\
4	0.6293	\\
5	0.6144	\\
6	0.6074	\\
7	0.5935	\\
8	0.5791	\\
9	0.5637	\\
10	0.5566	\\
11	0.5553	\\
12	0.5612	\\
13	0.585	\\
14	0.5993	\\
15	0.6225	\\
16	0.6359	\\
17	0.6694	\\
18	0.7137	\\
19	0.7534	\\
20	0.7791	\\
};

\addplot [color=red,dashed]
  table[row sep=crcr]{
12	0.00	\\
12	0.80	\\
};

\end{axis}
\end{tikzpicture}%

%% file: 2busgencost_20scen.tikz
%
%
%
\definecolor{mycolor1}{rgb}{0.00000,0.75000,0.75000}%
\definecolor{mycolor2}{rgb}{0.75000,0.00000,0.75000}%
\definecolor{mycolor3}{rgb}{0.75000,0.75000,0.00000}%
\begin{tikzpicture}

\begin{axis}[%
width=0.8\columnwidth,
height=2.2in,
scale only axis,
xmin=0,
xmax=50,
ymin=9,
ymax=18,
xlabel=\small Wind Penetration (\%),
ylabel=\small Cost of generation (\$MWh),
legend style={draw=black,fill=white,legend cell align=left}
]
\addplot [color=blue,solid, every mark/.append style={solid}, mark=*]
  table[row sep=crcr]{
0	17.63888	\\
5	17.133180861605	\\
10	16.62748172321	\\
15	16.121782584815	\\
20	15.6183075826982	\\
25	15.1404878117061	\\
30	14.662668040714	\\
35	14.1876970995525	\\
40	13.742609192701	\\
45	13.3386656211211	\\
50	12.947120443221	\\
};
\addlegendentry{ $\scriptsize \text{$\pm$0\%  Flexibility}$};

\addplot [color=black!50!green,dotted, every mark/.append style={solid}, mark=square*]
  table[row sep=crcr]{
0	17.63888	\\
5	16.8271591960815	\\
10	16.0154383921631	\\
15	15.2037175882446	\\
20	14.3919967843262	\\
25	13.5802759804077	\\
30	12.9097093740473	\\
35	12.4553600436134	\\
40	12.043953576799	\\
45	11.6524083988989	\\
50	11.2608632209988	\\
};
\addlegendentry{ $\scriptsize \text{$\pm$10\% Flexibility}$};

\addplot [color=red,densely dotted, every mark/.append style={solid}, mark=otimes*]
  table[row sep=crcr]{
0	17.63888	\\
5	16.8271591960815	\\
10	16.0154383921631	\\
15	15.2037175882446	\\
20	14.3919967843262	\\
25	13.5802759804077	\\
30	12.7685551764892	\\
35	11.9568343725708	\\
40	11.1451135686523	\\
45	10.3333927647339	\\
50	9.59414531315626	\\
};
\addlegendentry{ $\scriptsize \text{$\pm$20\% Flexibility}$};

\addplot [color=mycolor1,loosely dotted, every mark/.append style={solid}, mark=triangle*]
  table[row sep=crcr]{
0	17.63888	\\
5	16.8271591960815	\\
10	16.0154383921631	\\
15	15.2037175882446	\\
20	14.3919967843262	\\
25	13.5802759804077	\\
30	12.7685551764892	\\
35	11.9568343725708	\\
40	11.1451135686523	\\
45	10.3333927647338	\\
50	9.52167196081539	\\
};
\addlegendentry{ $\scriptsize \text{$\pm$30\% Flexibility}$};



\end{axis}
\end{tikzpicture}%

%% file: windspill.tikz
%
%

\begin{tikzpicture}

\begin{axis}[%
width=0.8\columnwidth,
height=1.5in,
scale only axis,
xmin=0,
xmax=50,
ymin=2,
ymax=11,
xlabel=Load Flexibility (\%),
ylabel=Mean Wind spillage (MW),
xtick={0 ,10 ,20 ,30 ,40 ,50},
legend style={draw=black,fill=white,legend cell align=left}
]
\addplot [color=blue,dashed,mark=asterisk,mark options={solid}]
  table[row sep=crcr]{
0	8.53318527577578	\\
10	4.25425072845001	\\
20	2.78533728080236	\\
30	2.76180697160079	\\
40	2.76180697160079	\\
50	2.76180697160079	\\
};
\addlegendentry{$\scriptsize \text{20  Scenarios}$};

\addplot [color=red,dashed,mark=asterisk,mark options={solid}]
  table[row sep=crcr]{
0	10.152746804632	\\
10	5.56258864566067	\\
20	3.74390357597984	\\
30	3.67793571140596	\\
40	3.67793571140596	\\
50	3.67793571140597	\\
};
\addlegendentry{$\scriptsize \text{100 Scenarios}$};

\end{axis}
\end{tikzpicture}%

%% file: valueofstochasticsol.tikz
%
%
%
\definecolor{mycolor1}{rgb}{0.00000,0.75000,0.75000}%
\definecolor{mycolor2}{rgb}{0.75000,0.00000,0.75000}%
\definecolor{mycolor3}{rgb}{0.75000,0.75000,0.00000}%
\begin{tikzpicture}

\pgfplotscreateplotcyclelist{my black white}{%
solid, every mark/.append style={solid, fill=gray}, mark=*\\%
dotted, every mark/.append style={solid, fill=gray}, mark=square*\\%
densely dotted, every mark/.append style={solid, fill=gray}, mark=otimes*\\%
loosely dotted, every mark/.append style={solid, fill=gray}, mark=triangle*\\%
dashed, every mark/.append style={solid, fill=gray},mark=diamond*\\%
loosely dashed, every mark/.append style={solid, fill=gray},mark=*\\%
densely dashed, every mark/.append style={solid, fill=gray},mark=square*\\%
dashdotted, every mark/.append style={solid, fill=gray},mark=otimes*\\%
dasdotdotted, every mark/.append style={solid},mark=star\\%
densely dashdotted,every mark/.append style={solid, fill=gray},mark=diamond*\\%
}

\begin{axis}[%
width=0.8\columnwidth,
height=1.5in,
scale only axis,
xmin=0,
xmax=50,
ymin=0,
ymax=9.0,
xlabel=Wind Penetration (\%),
ylabel=Difference in generation costs (\%),
legend style={draw=black,fill=white,legend cell align=left},
legend pos=north west
]
\addplot [color=blue,solid, every mark/.append style={solid}, mark=*]
  table[row sep=crcr]{
0	0	\\
5	0.60096849086074	\\
10	1.23849203594953	\\
15	1.91601055129166	\\
20	2.62279338138651	\\
25	3.2014925453489	\\
30	3.81790843333819	\\
35	4.45487422286181	\\
40	4.9681614246861	\\
45	5.25048725736949	\\
50	5.44936982889317	\\
};
\addlegendentry{ $\scriptsize \text{$\pm$0\%  Flexibility}$};

\addplot [color=black!50!green,dotted, every mark/.append style={solid}, mark=square*]
  table[row sep=crcr]{
0	2.06247721345783e-14	\\
5	0.283747173158243	\\
10	0.596257028656674	\\
15	0.94213645309939	\\
20	1.32719282736988	\\
25	2.19445246330884	\\
30	4.25166741126573	\\
35	4.85524676758253	\\
40	5.22910759502616	\\
45	5.44936982889281	\\
50	5.68494929268361	\\
};
\addlegendentry{ $\scriptsize \text{$\pm$10\% Flexibility}$};

\addplot [color=red,densely dotted, every mark/.append style={solid}, mark=otimes*]
  table[row sep=crcr]{
0	2.06247721345783e-14	\\
5	0.283747173158264	\\
10	0.5962570286568	\\
15	0.942136453099509	\\
20	1.32703166231197	\\
25	1.7579388154782	\\
30	2.24363310719701	\\
35	2.79527263971163	\\
40	3.42768241633985	\\
45	4.5699962456554	\\
50	5.89537455647018	\\
};
\addlegendentry{ $\scriptsize \text{$\pm$20\% Flexibility}$};

\addplot [color=mycolor1,loosely dotted, every mark/.append style={solid}, mark=triangle*]
  table[row sep=crcr]{
0	2.06247721345783e-14	\\
5	0.283747173158329	\\
10	0.596257028656823	\\
15	0.942136453099426	\\
20	1.32703166231197	\\
25	1.75793881547832	\\
30	2.24363310719711	\\
35	2.79527263971164	\\
40	3.42726618244839	\\
45	4.15854992314018	\\
50	5.01451727577147	\\
};
\addlegendentry{ $\scriptsize \text{$\pm$30\% Flexibility}$};



\end{axis}
\end{tikzpicture}%

%% file: LMPS_39bus.tikz
%
%
%
\definecolor{mycolor1}{RGB}{219,144,71}
\definecolor{mycolor2}{RGB}{186,146,162}

\begin{tikzpicture}

\begin{axis}[%
width=0.8\columnwidth,
height=1.8in,
scale only axis,
xmin=0,
xmax=12,
ymin=1540,
ymax=1860,
xlabel={$t$},
ylabel={\$/MWh},
legend style={draw=black,fill=white,legend cell align=left}
]
\addplot [color=blue,solid, every mark/.append style={solid}, mark=*]
  table[row sep=crcr]{
1	1730.02921679187	\\
2	1714.28857922213	\\
3	1697.62991112124	\\
4	1671.78466750401	\\
5	1640.62218493177	\\
6	1621.91948474411	\\
7	1615.84897510769	\\
8	1601.11344103346	\\
9	1576.72648138402	\\
10	1572.13943967613	\\
11	1566.61514768961	\\
12	1587.76951346926	\\
};
\addlegendentry{$\scriptsize \text{$\pm$0\%  Load Flexibility}$};

\addplot [color=black!50!green,dotted, every mark/.append style={solid}, mark=square*]
  table[row sep=crcr]{
1	1708.05649667839	\\
2	1692.30591249592	\\
3	1676.07064257164	\\
4	1650.9664912268	\\
5	1620.29680266531	\\
6	1601.82579471799	\\
7	1596.21508370891	\\
8	1581.95590018853	\\
9	1558.0783884477	\\
10	1553.72622821918	\\
11	1548.24522734578	\\
12	1569.20413456889	\\
};
\addlegendentry{$\scriptsize \text{$\pm$10\% Load Flexibility}$};

\addplot [color=red,densely dotted, every mark/.append style={solid}, mark=otimes*]
  table[row sep=crcr]{
1	1590.96729914139	\\
2	1590.96729914378	\\
3	1590.96729914416	\\
4	1590.96729913632	\\
5	1590.96729913618	\\
6	1590.96729915179	\\
7	1590.96729914258	\\
8	1590.96729913909	\\
9	1590.96729896244	\\
10	1565.36133864638	\\
11	1562.20183581729	\\
12	1576.65174348595	\\
};
\addlegendentry{$\scriptsize \text{$\pm$20\% Load Flexibility}$};

\addplot [color=purple,loosely dotted, every mark/.append style={solid}, mark=triangle*]
  table[row sep=crcr]{
1	1598.96500481969	\\
2	1598.96500482329	\\
3	1598.96500482506	\\
4	1598.9650048123	\\
5	1598.96500481142	\\
6	1598.96500481207	\\
7	1598.96500481051	\\
8	1598.96500480863	\\
9	1598.96500479153	\\
10	1583.7745539511	\\
11	1580.57205841721	\\
12	1595.21711691523	\\
};
\addlegendentry{$\scriptsize \text{$\pm$30\% Load Flexibility}$};

\addplot [color=mycolor2,dashed, every mark/.append style={solid},mark=diamond*]
  table[row sep=crcr]{
1	1604.01477365026	\\
2	1604.01477365104	\\
3	1604.01477365139	\\
4	1604.01477364979	\\
5	1604.01477364958	\\
6	1604.01477364955	\\
7	1604.01477364941	\\
8	1604.01477364926	\\
9	1604.01477364869	\\
10	1602.18775898043	\\
11	1598.94223141653	\\
12	1604.01477364679	\\
};
\addlegendentry{$\scriptsize \text{$\pm$40\% Load Flexibility}$};

\addplot [color=black,loosely dashed, every mark/.append style={solid},mark=star]
  table[row sep=crcr]{
1	1609.09350402844	\\
2	1609.09350401698	\\
3	1609.09350400826	\\
4	1609.09350399188	\\
5	1609.09350398549	\\
6	1609.09350398537	\\
7	1609.09350398086	\\
8	1609.09350397522	\\
9	1609.09350396235	\\
10	1609.0935039447	\\
11	1609.0935039404	\\
12	1609.09350395226	\\
};
\addlegendentry{$\scriptsize \text{$\pm$100\% Load Flexibility}$};
\end{axis}
\end{tikzpicture}%

%% file: loadprofile_39bus.tikz
%
%
%
\definecolor{mycolor1}{RGB}{219,144,71}
\definecolor{mycolor2}{RGB}{186,146,162}
\begin{tikzpicture}

\begin{axis}[%
width=0.8\columnwidth,
height=1.8in,
scale only axis,
xmin=0,
xmax=12,
ymin=0.54,
ymax=0.68,
xlabel={$t$},
ylabel={Demand (p.u.)},
legend style={draw=black,fill=white,legend cell align=left}
]
\addplot [color=blue,solid, every mark/.append style={solid}, mark=*]
  table[row sep=crcr]{
1	0.663899470092737	\\
2	0.664199334352035	\\
3	0.651405125955186	\\
4	0.629015261260701	\\
5	0.614122003048744	\\
6	0.607125170331717	\\
7	0.593231459650766	\\
8	0.578837975204311	\\
9	0.563444943226852	\\
10	0.556348155756724	\\
11	0.55504874396642	\\
12	0.560946074399341	\\
};

\addplot [color=black!50!green,dotted, every mark/.append style={solid}, mark=square*]
  table[row sep=crcr]{
1	0.662867790711217	\\
2	0.661363420433768	\\
3	0.651400640152109	\\
4	0.631258047104943	\\
5	0.614044387042293	\\
6	0.604662254420617	\\
7	0.593825102086112	\\
8	0.580417626229313	\\
9	0.56353602540237	\\
10	0.557287800975636	\\
11	0.554900154159611	\\
12	0.562060468543974	\\
};

\addplot [color=red,densely dotted, every mark/.append style={solid}, mark=otimes*]
  table[row sep=crcr]{
1	0.650168418598996	\\
2	0.647135222782816	\\
3	0.640781316782394	\\
4	0.625335405001929	\\
5	0.605230222222547	\\
6	0.598736798020263	\\
7	0.59563964759072	\\
8	0.590350583969208	\\
9	0.57986462069983	\\
10	0.567663001003888	\\
11	0.565439716557886	\\
12	0.571278764165854	\\
};

\addplot [color=purple,loosely dotted, every mark/.append style={solid}, mark=triangle*]
  table[row sep=crcr]{
1	0.633697902359506	\\
2	0.630523257546618	\\
3	0.62494041634237	\\
4	0.610747467205555	\\
5	0.603697935763436	\\
6	0.603564972413707	\\
7	0.601942886280543	\\
8	0.598345147336754	\\
9	0.589456604095029	\\
10	0.579690281430829	\\
11	0.577475795645022	\\
12	0.58354105101553	\\
};

\addplot [color=mycolor2,dashed, every mark/.append style={solid},mark=diamond*]
  table[row sep=crcr]{
1	0.616417527571899	\\
2	0.613219487054141	\\
3	0.608096777668893	\\
4	0.605684507112426	\\
5	0.606729853292923	\\
6	0.606992367509512	\\
7	0.606207128974961	\\
8	0.603973082487244	\\
9	0.597039655403093	\\
10	0.591719426989833	\\
11	0.589404834031134	\\
12	0.592139069162852	\\
};

\addplot [color=black,loosely dashed, every mark/.append style={solid},mark=star]
  table[row sep=crcr]{
1	0.585237716187722	\\
2	0.593206722828717	\\
3	0.601406877673971	\\
4	0.606732672226076	\\
5	0.608435850435781	\\
6	0.608892630230534	\\
7	0.608832001005479	\\
8	0.608446714841414	\\
9	0.606147896967919	\\
10	0.605283461166792	\\
11	0.604155638084479	\\
12	0.600845535731606	\\
};
\end{axis}
\end{tikzpicture}%

%% file: times.tikz
%
%

\begin{tikzpicture}

\begin{axis}[%
width=0.8\columnwidth,
height=2.05in,
scale only axis,
xmin=0.5,
xmax=8.5,
ymode=log,
ymin=0.01,
ymax=100,
legend style={draw=black,fill=white,legend cell align=left},
legend pos=north west,
ylabel={$t~(s)$},
xtick={1,2,3,4,5,6,7,8},
xticklabels={9,14,24,30,39,57,118,300},
xlabel={$n^\mc{B}$}
]
\addplot [color=orange,only marks,mark=asterisk,mark size=2pt]
  table[row sep=crcr]{
1	0.04	\\
2	0.10	\\
3	0.78	\\
4	0.18	\\
5	0.24	\\
6	1.6	\\
7	7.86	\\
8	14.0	\\
};
\addplot [color=blue,only marks,mark=o,mark options={solid},mark size=1.5pt]
  table[row sep=crcr]{
1	0.054	\\
2	0.154	\\
3	1.01	\\
4	0.2358	\\
5	0.282	\\
6	1.7685	\\
7	9.1956	\\
8	15.5191	\\
};

\addplot [color=red,only marks, mark=triangle*,mark size=2pt]
  table[row sep=crcr]{
1	0.06	\\
2	0.18	\\
3	1.34	\\
4	0.316	\\
5	0.32	\\
6	2.084	\\
7	10.8	\\
8	19.73	\\
};
\end{axis}
\end{tikzpicture}%